\documentclass[11pt]{article}
\usepackage{arxiv}
\usepackage{cite}
\usepackage{amsmath,amssymb,amsfonts}
\usepackage{algorithmic}
\usepackage{graphicx}
\usepackage{textcomp}
\usepackage{xcolor}
\usepackage[capitalise]{cleveref}
\usepackage{hhline}

\def\BibTeX{{\rm B\kern-.05em{\sc i\kern-.025em b}\kern-.08em
    T\kern-.1667em\lower.7ex\hbox{E}\kern-.125emX}}
\begin{document}

\title{Fast Solving Complete 2000-Node Optimization Using Stochastic-Computing Simulated Annealing\\
%{\footnotesize \textsuperscript}
}

\author{
  Kota Katsuki$^{1,2}$,
  Duckgyu Shin$^{1,2}$,
  Naoya Onizawa$^{2}$,
  Takahiro Hanyu$^{2}$ \\
  \texttt{\{kota.katsuki.t5@tohoku.ac.jp, duckgyu.shin.p4@dc.tohoku.ac.jp,} \\
  \texttt{naoya.onizawa.a7@tohoku.ac.jp, takahiro.hanyu.c4@tohoku.ac.jp\}} \\
  $^{1}$ Graduate School of Engineering, Tohoku University, Sendai, Japan \\
  $^{2}$ Research Institute of Electrical Communication, Tohoku University, Sendai, Japan
}
\date{}

%\author{\IEEEauthorblockN{Kota Katsuki}
%\IEEEauthorblockA{\textit{Tohoku University} \\
%\textit{Graduate School of Engineering}\\
%Miyagi, Japan \\
%kota.katsuki.t5@tohoku.ac.jp}
%\and
%\IEEEauthorblockN{Duckgyu Shin}
%\IEEEauthorblockA{\textit{Tohoku University} \\
%	\textit{Graduate School of Engineering}\\
%	Miyagi, Japan \\
%	duckgyu.shin.p4@dc.tohoku.ac.jp}
%\and
%\IEEEauthorblockN{Naoya Onizawa}
%\IEEEauthorblockA{\textit{Tohoku University} \\
%\textit{Research Institute for Electrical Communication}\\
%Miyagi, Japan \\
%naoya.onizawa.a7@tohoku.ac.jp}
%\and
%\IEEEauthorblockN{Takahiro Hanyu}
%\IEEEauthorblockA{\textit{Tohoku University} \\
%\textit{Research Institute for Electrical Communication}\\
%Miyagi, Japan \\
%takahiro.hanyu.c4@tohoku.ac.jp}
%}

\maketitle

\begin{abstract}
In this paper, we evaluate stochastic-computing simulated annealing (SC-SA) for solving large-scale combinatorial optimization problems.
SC-SA is designed using stochastic computing, where the computatoin is reazlied using random bitstream, resulting in fast converging to the global minimum energy of the problems.
The proposed SC-SA is compared with a typical SA and existing simulated-annealing (SA) processors on the maximum cut (MAX-CUT) problems, such as Gset that is a benchmark for SA.
The simulation results show that SC-SA realizes a few orders of magnitude faster than a typical SA.
In addition, SC-SA achieves better MAX-CUT scores than other existing methods on K2000 that is a complete 2000-node optimization problem.
\end{abstract}


\section{Introduction}
\label{sec:introduction}
Simulated annealing (SA) is a well-known possible technique for solving NP-hard problems \cite{NP-hard}.
Various SA techniques have been  studied  for high-speed solvers of Ising models on combinatorial optimization problems  \cite{Ising_HW1,Ising_HW2}.
Among them, stochastic computing-simulated annealing (SC-SA) has been recently presented for a fast converging SA algorithm \cite{Onizawa}, 
it can solve the combinatorial optimizatoin problems using probabilistic bits (p-bits) that are approximated by stochastic computing, where the p-bits were originally presented for  probabilistic bidirectional computing between inputs and outputs \cite{IL}.

SC-SA is designed based on  stochastic computing and an extension of stochastic computing, \emph{integral stochastic computing} \cite{SDNN}.
Stochastic computing  is a probabilistic computing technique using random bit streams \cite{stochastic_first,stochastic}, and integral stochastic computing is a computing system that expresses a numerical value by the probability of existence of ``1'' in the bit stream.
Stochastic computing has the advantage that it can be realized by approximating a complicated operation circuit, such as the $\tanh$ function with a CMOS circuit.
The $\tanh$ function in the p-bit function is approximated using a saturated up-down counter that can induce relaxed transitions of a spin state for fast convergence to the global minimum energy.
Its effectiveness was confirmed in comparison with quantum annealing (D-wave) and a conventional SA algorithm \cite{Onizawa}.
On the other hand, the combinatorial optimization problem evaluated above is small, and it has not been verified to be effective on a large scale.

In this paper, we evaluate the proposed SC-SA method on a large-scale maximum cut (MAX-CUT) problem, which contains 2000-vertex at most, and compare with a conventional SA and existing SA processors.
The MAX-CUT problem is a typical combinatorial optimization problem as a large-scale problem.
The proposed SC-SA method finds the near-optimal solution and has better quality than that of the conventional SA on Gset and K2000 that are benchmarks on MAX-CUT for SA.
The simulation results show that the proposed SC-SA achieves a convergence speed that is a few orders of magnitude faster.
In addition, we compare the quality of the solution found by the proposed SC-SA method with the near-optiomal solution from the existing annealing processors.
This paper is organized as follows:
\cref{sec:SA} introduces simulated annealing based on integral stochastic computing.
\cref{sec:setup} describes the simulation conditions of SA for combinatorial optimization problems.
\cref{sec:evaluation} compares the proposed SC-SA method with conventional SA and existing annealing processor, and discusses some insights into the proposed simulated annealing method.
\cref{sec:conclusion} concludes the brief.
\section{Simulated Annealing Based on Integral Stochastic Computing}
\label{sec:SA}
\begin{figure}[tb]
	\begin{center}
		\includegraphics[width=0.8\linewidth]{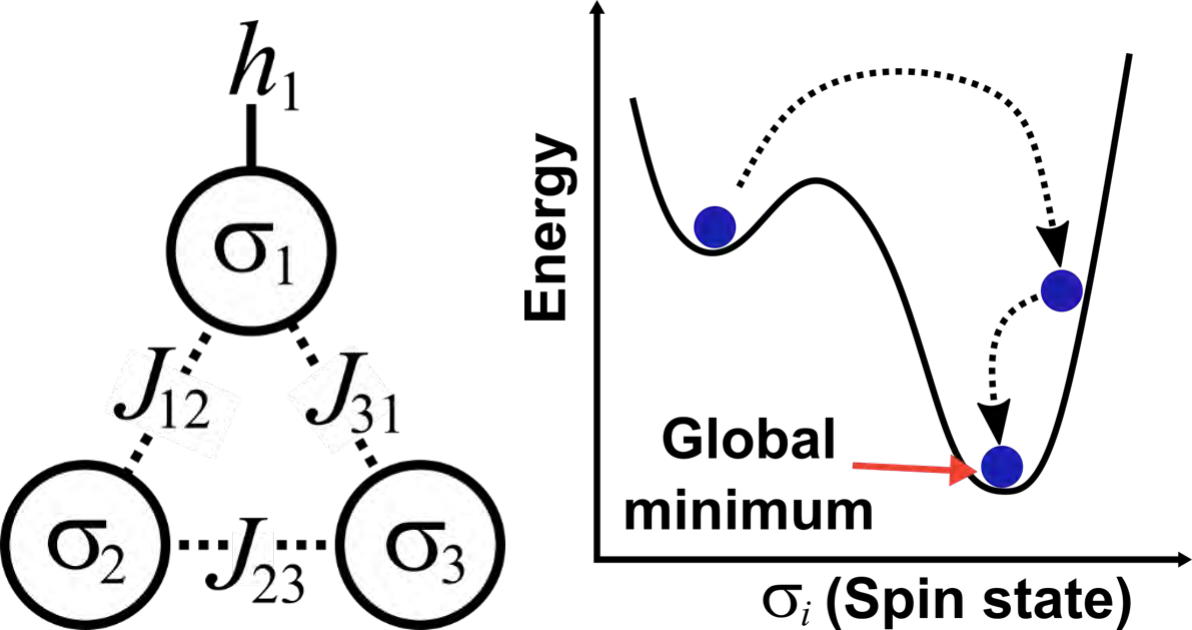}
		\caption{Example of an Ising model with three spins. The Ising model is used to represent combinatorial optimization problems solved by SA or SC-SA with reaching the global minimum energy of the model.}
		\label{fig:f1}
	\end{center}
\end{figure}
The SC-SA algorithm is the SA method on the Ising model which is converted from the combinatorial optimization problem \cite{Ising}.
The Ising model can be represented as a network of the spin, and \cref{fig:f1} shows the Ising model with 3-spin.
SC-SA searches a solution of the problem by converging Hamiltonian ($H$) which is an energy of the Ising model, and $H$ is given by:
\begin{equation}
	H = - \sum_i \sigma_i h_i - \frac{1}{2} \sum_{i \neq j} \sigma_i \sigma_j J_{ij},
\label{eqn:energy}
\end{equation}
where $h_i$ is a bias of the spin, $J_{ij}$ is a weight of the interconnection between the spins, and $\sigma_i \in [-1, +1]$ is the state of the spin.
$h_i$ and $J_{ij}$ are determined during converting the optimization problem to the Ising model, so that the optimal solution of the optimization problem is embedded at the global minimum of the Hamiltonian.

While the annealing process of SC-SA, $H$ is converged to the global minimum by fluctuating the state of the spins using random noise signals.
The spins are changes its state probabilistically, and its calculations are defined as:
\begin{subequations}
	\begin{equation}
		I_i(t+1) = h_i + \sum_j J_{ij} \cdot \sigma_j (t) + n_{rnd} \cdot r_i (t),
		\label{eqn:I}
	\end{equation}
	\begin{equation}
		\mathrm{Itanh}_i (t+1) =
		\begin{cases}
			I_0-1, \text{if} \ \mathrm{Itanh}_i (t) + I_i (t+1) \geq I_0 \\
			-I_0, \text{else if} \ \mathrm{Itanh}_i (t) + I_i (t+1) < -I_0 \\
			\mathrm{Itanh}_i (t) + I_i (t+1), \ \text{otherwise}.
		\end{cases}
		\label{eqn:updown}
	\end{equation}
	\begin{equation}
		\sigma_i(t+1)= \mathrm{sgn}(\mathrm{Itanh}_i(t+1)) = 
		\begin{cases}
			1, & \text{if} \ \mathrm{Itanh}_i(t+1) \geq 0 \\
			-1, & \text{otherwise},
		\end{cases}
		\label{eqn:m}
	\end{equation}
	\label{eqn:prop}
\end{subequations}
\begin{figure}[t]
	\centering
	\includegraphics[width=0.8\linewidth]{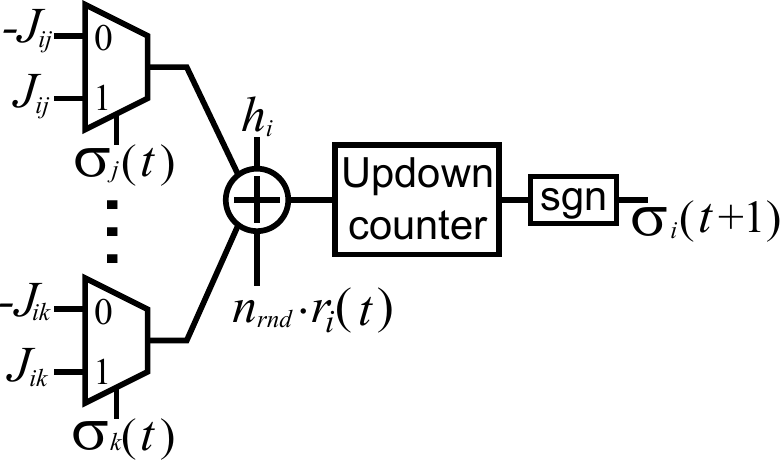}
	\caption{Spin-gate circuit that calculates \cref{eqn:I}, \cref{eqn:updown}, and \cref{eqn:m} based on stochastic computing and integral stochastic computing.}
	\label{fig:spin_gate}
\end{figure}
where $t$ is a cycle, $n_{rnd}$ is a magnitude of the noise signals, $r_i(t)$ is the random noise signals, $\mathrm{sgn}$ is a sign function, and $I_0$ is a pseudoinverse temperature \cite{Onizawa}.
In SC-SA, \cref{eqn:prop} is realized using a spin-gate circuit shown in \cref{fig:spin_gate}.
SC-SA is based on stochastic computing, thus calculations of the $\mathrm{Itanh}$ can be realized using the saturated up-down counter \cite{stochastic_first, stochastic_book}.
Also, multiplications of $\sigma_i$ and $J_{ij}$ is performed by a multiplexer, and additions is performed by a binary adder.
The pseudoinverse temperature. $I_0$ controls the converting speed of the Hamiltonian.
For obtaining the solution which has high quality, $I_0$ need to be controlled precisely.
The detail of the pseudoinverse temperature is introduced in \cref{sec:setup}.
\section{Simulation Setup}
\label{sec:setup}
\subsection{Hamiltonian design of MAX-CUT Problem}
In this paper, we apply SC-SA to solve the MAX-CUT problem that is a typical combinatorial optimization problem.
The MAX-CUT problem is defined on a graph $G=(V, E)$, where $V = \{1,2,...,n\}$ and $E \subset \{(i,j):1\leq i<j\leq n\}$ are given.
Let the edge weights $w_{ij} = w_{ji}$ be given such that $w_{ij} = 0$ for $(i,j) \notin E$.
The MAX-CUT problem aims to find a bipartition ($V_1$, $V_2$) of $V$ so that the sum of the edge weights between $V_1$ and $V_2$ is maximized. 
From \cite{QUBO_maxcut}, the argument $h_i$ of the Ising model of the MAXCUT problem with $n$ nodes can be represented by a 0 matrix of $(n \times 1)$, and $J_{ij}$ can be represented by an adjacent matrix of $(n \times n)$ of the input graph.
\begin{figure}[t]
	\centering
	\includegraphics[width=0.7\linewidth]{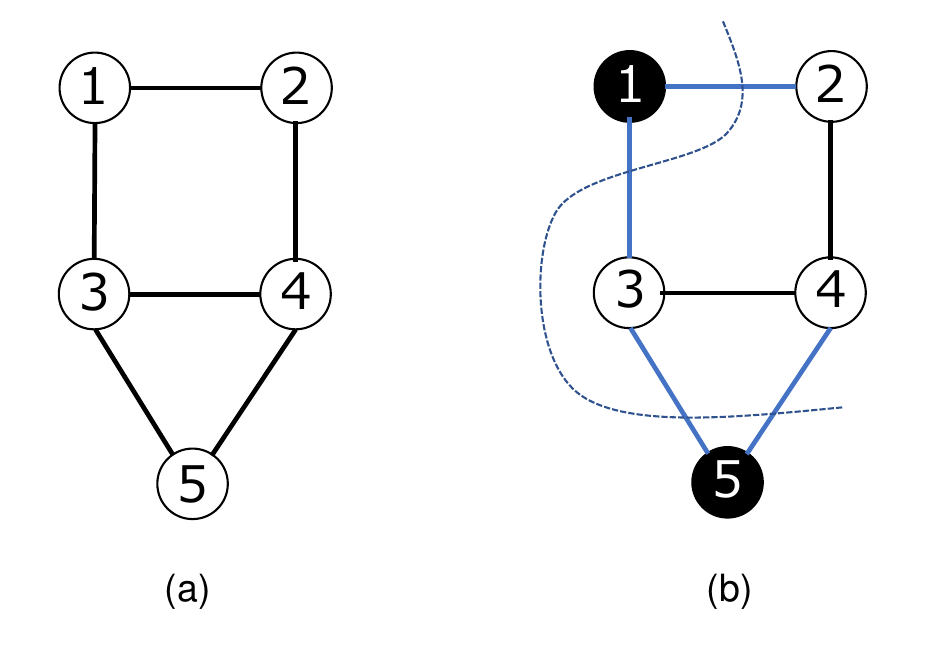}
	\caption{(a) Example of a MAX-CUT problem with five veritces and (b) the solution of MAX-CUT problem. A cut value of the solution is maximized with two separated groups of Group A (veritces 1,5) and Group B (veritces 2,3,4).}
	\label{fig:problem}
\end{figure}
\begin{figure}[t]
	\centering
	\includegraphics[width=0.9\linewidth]{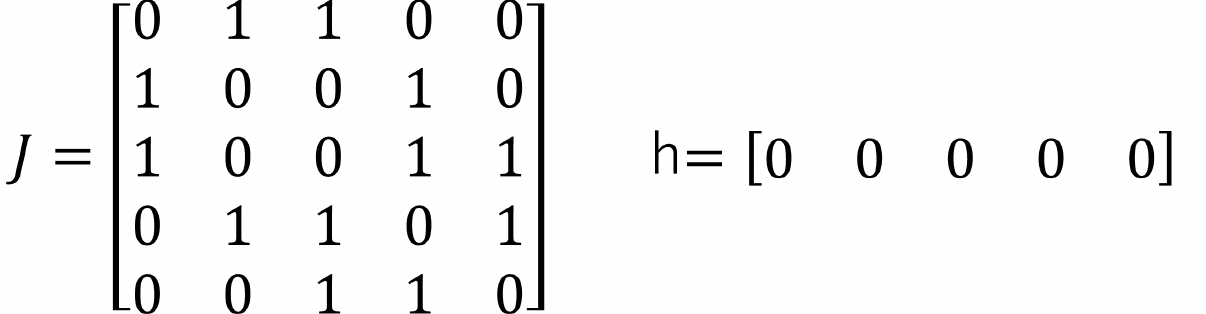}
	\caption{$h$ and $J$ for the MAX-CUT problem shown in in Fig.3 (a).}
	\label{fig:example}
\end{figure}

\cref{fig:problem} shows an example of a five-veritce MAX-CUT problem with weights of +1.
The graph is divided into Group A (veritces 1,5) and Group B (veritces 2,3,4), where the sum of the edge weights is 4.
Fig. 4 shows adjacent matrices $J$ and $h$ of the MAX-CUT problem shown in fig. 3 (a).
\begin{figure}[t]
	\centering
	\includegraphics[width=0.8\linewidth]{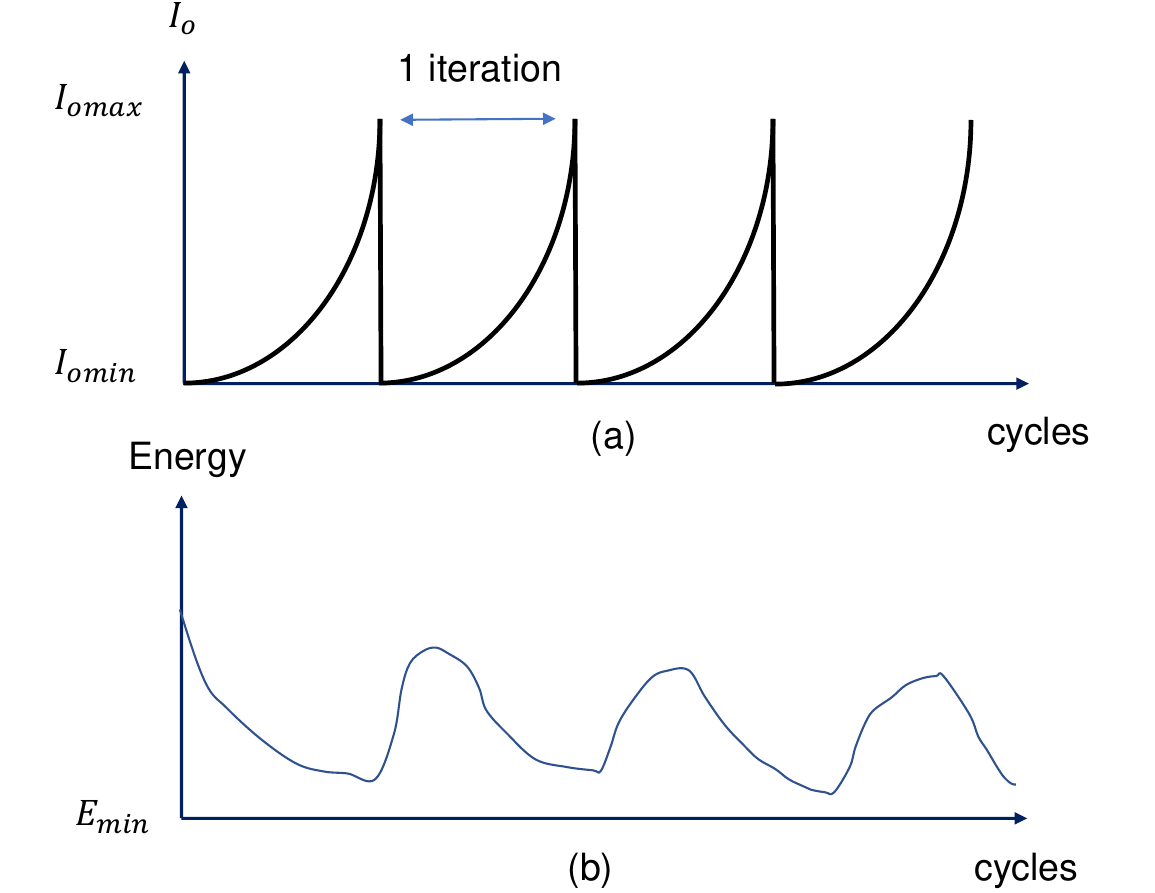}
	\caption{Simulation conditions: (a) $I_0$ control as $I_0(t+\tau) = (1/\beta)\cdot I_0(t)$ for reaching the global minimum energy and (b) an example of an energy transition.}
	\label{fig:condition}
\end{figure}
\subsection{Simulation condition}
In the proposed SC-SA method, \cref{eqn:prop} are calculated to obtain the global minimum energy of the Ising model.
To reach the global minimum energy, a pseudoinverse temperature, $I_0$, is controlled within the range $I_{0min}$ to $I_{0max}$ for each iteration shown in \cref{fig:condition}.
In each iteration, $I_0$ is gradually increased according to $I_0(t+\tau) = (1/\beta)\cdot I_0(t)$.
When $I_0$ is small, the spin state, $\sigma_i$, is easily flipped between `-1' and `+1'.
In contrast, it become hard to flipped, thus, the state of nodes is stabilized.

The conventional SA process is as follows:
The temperature $T$ gradually decreases with each cycle depending on the $\Delta T$ as $T \gets 1 / (1 / T + \Delta T)$, where the initial temperature is 10,000 and the final temperature is 1, and the number of cycles varies with the $\Delta T$.
At each cycle, the state of the vertex is randomly inverted, and if the new energy ($E_{new}$) is lower or higher than the current energy ($E_{cur}$), the new state is accepted with a probability of $\exp(-(E_{new} - E_{cur}) / T)$.

The conventional and proposed SA methods are evaluated using MATLAB R2020b on an 8-core Intel Core i7 at 2.3 GHz and 32 GB memory.

\section{Evaluation}
\label{sec:evaluation}
\subsection{Gset}
SA and SC-SA are applied to solve data sets of MAX-CUT problems, \emph{Gset} \cite{Gset}.
G6, G11, and G16 in Gset have 800 vertices as shown in the \cref{tb:Gset_spec}; however, the edges between vertices and the weights of edges are different.
\begin{table}[t]
	\caption{Spectification of the MAX-CUT problems used for simulations.}
	\label{tb:Gset_spec}
	\centering
	\begin{tabular}{|c||c|c|c|c|}
		\hline
		& Node & Edge & Weight & Graph type \\
		\hhline{|=#=|=|=|=|}
		G6 & 800 & 19,176 & \{-1, +1\} &Ramdom\\
		\hline
		G14 & 800 & 4,694 & \{+1\} &Troidal\\
		\hline
		G18 & 800 & 4,694 & \{-1, +1\} &Planar\\
		\hline
		K2000 & 2,000 & 1,999,000 & \{-1, +1\} &Ramdom\\
		\hline
	\end{tabular}
\end{table}
In order to obtain the global minimum energy with high probability in SC-SA, a parameter search is performed for $\tau$, $I_{0min}$ and $n_{rnd}$ while $I_{0min}$ and $\beta$ are fixed.
\cref{tb:parameter} shows the results of parameter search for G6, G14, and G16.
\begin{table}[t]
	\caption{Summary of the simulation parameters for SC-SA.}
	\label{tb:parameter}
	\centering
	\begin{tabular}{|c||c|c|c|}
		\hline
		 & $\tau$  & $I_{0max}$ & $n_{rnd}$ \\
		\hhline{|=#=|=|=|}
		G6    & 500 & 64 & 4\\
		\hline
		G14   & 1 & 512 & 2\\
		\hline
		G18   & 500 & 32 & 4\\
		\hline
		K2000 & 500 & 1024 & 32\\
		\hline
	\end{tabular}
\end{table}

%\cref{tb:Gset} shows the average edge weight vs the number of cycles for 100 trials.

\cref{tb:Gset} shows the average edge weight at 1,000 cycles for 100 trials.
SC-SA scored better than conventional SA for any data sets.

\begin{table}[t]
	\caption{Comparison of annealing  accuracy on Gset}
	\label{tb:Gset}
	\centering
	\begin{tabular}{|c||c|c|}
		\hline
		& SA & SC-SA \\
		\hhline{|=#=|=|}
		G6 & 1,935 & 2003.5 \\
		\hline
		G14 & 2,801 & 2,928.7 \\
		\hline
		G18 & 886.53 & 932.4 \\
		\hline
	\end{tabular}
\end{table}

\subsection{K2000}
\begin{figure}[t]
	\centering
	\includegraphics[width=0.8\linewidth]{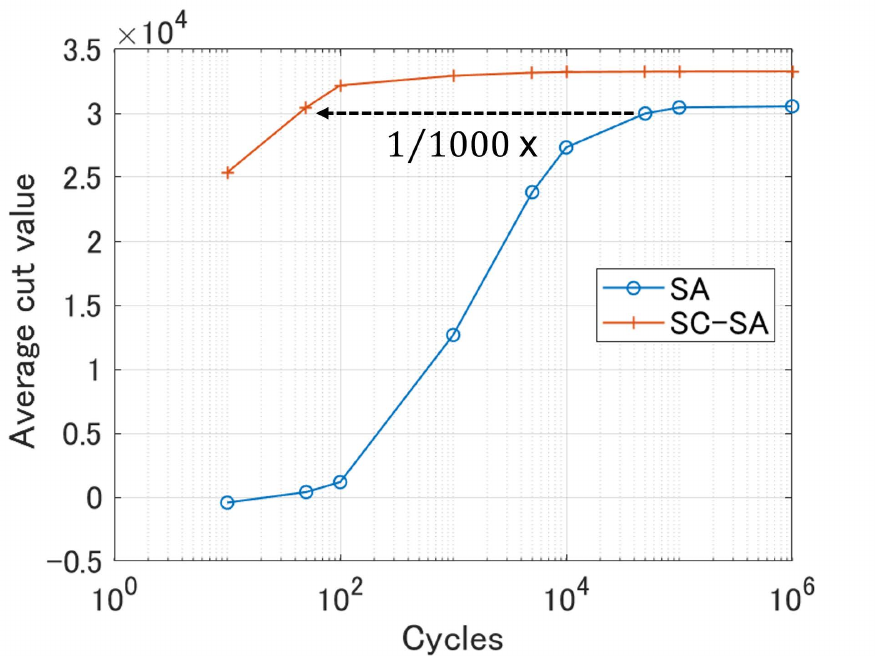}
	\caption{The average cut value vs. the number of cycles on the MAX-CUT problem K2000. The performance is evaluated with 100 trials each for the conventional SA and the proposed SC-SA methods.}
	\label{fig:cycle_maxcut}
	\vspace{-3mm}
\end{figure}

\begin{figure}[t]
	\centering
	\includegraphics[width=0.8\linewidth]{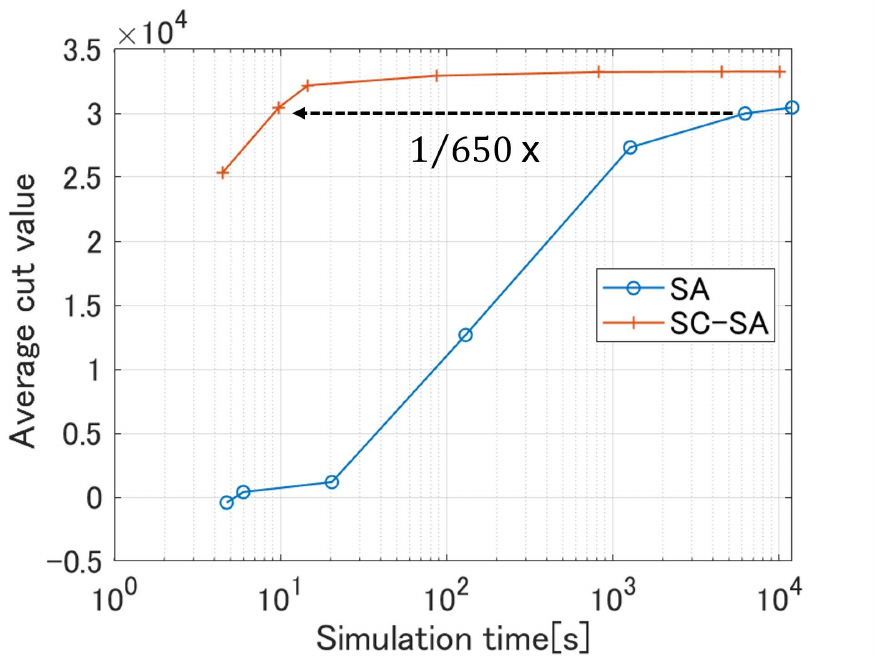}
	\caption{The average cut value vs. the simulation time on the MAX-CUT problem K2000. The performance is evaluated with 100 trials each for the conventional SA and the proposed SC-SA methods.}
	\label{fig:speed_maxcut}
	\vspace{-3mm}
\end{figure}

In this subsection, the MAXCUT problem (K2000) is solved.
K2000 is a complete graph which contains 2,000 vertices, and its number of edges is 1,999,000 \cite{CIM}.
% The K2000 contains 2000 vertices and 1999000 edges.
Edge weights can be either -1 or + 1.
The result of parameter search for the K2000 is shown in the \cref{tb:parameter}.
\cref{fig:cycle_maxcut} shows the average cut value vs. the number of cycle for 100 trials.
As a result, the conventional SA takes 50,000 cycles to obtain cut value as 30,000, which is approximately 90\% of the near-optimal solution (33337) \cite{best_K2000}.
Next, the average cut value vs. the simulation time is evaluated in \cref{fig:speed_maxcut}.
From \cref{fig:speed_maxcut}, SC-SA is approximately 650 times faster than the conventional SA in obtaining the near-optimal solution.
Fig.8 shows the energy vs.  cycles , where the energy is the Hamiltonian defined in \cref{eqn:energy}.
The energy of SC-SA is rapidly dropped to reach the global minimum energy,  which can search for solutions with lower energy than the conventional SA method.

\begin{figure}[t]
	\centering
	\includegraphics[width=0.8\linewidth]{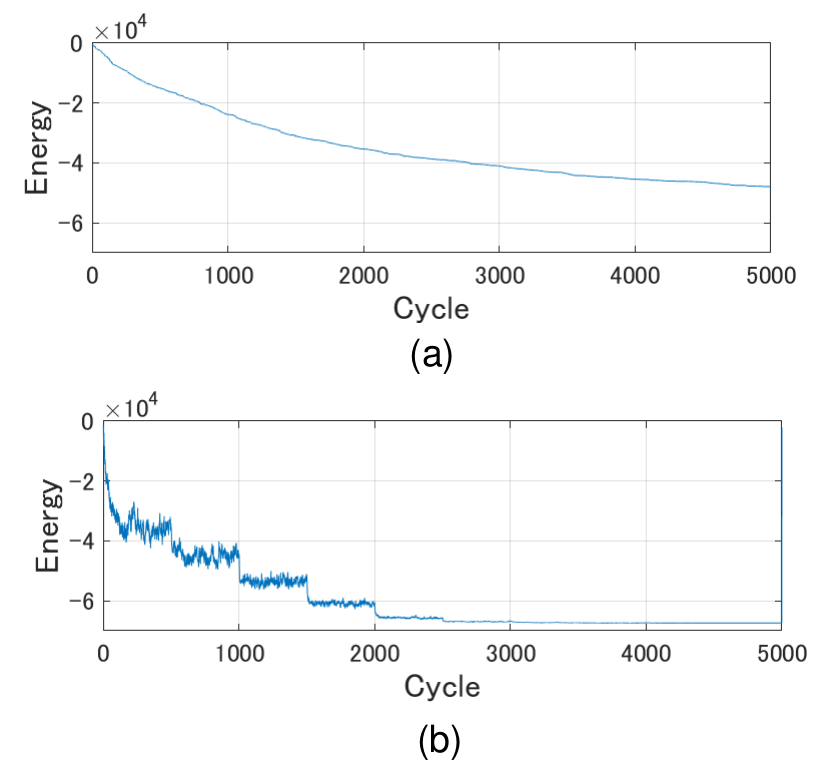}
	\caption{Energy vs. cycle  on simulating K2000 (a) the conventional SA (b) SC-SA. The energy of SC-SA is rapidly dropped to reach the global minimum energy.}
	\label{fig:Evscycles}
	\vspace{-3mm}
\end{figure}

The average cut value and the maximum cut value are shown in \cref{tb:annealing processors} in comparisons with existing annealing processors such as CIM \cite{CIM}, SB \cite{SB} and STATICA \cite{STATICA}.
The results of the proposed SC-SA method is the average cut value and the maximum cut value at 100,000 cycles.
Compared with existing annealing processors, the proposed SC-SA method achieves the best average cut value.
Compared with SB, which has the best average cut value among the three existing annealing processors, SC-SA obtained a score of 191 better than SB.
Though the existing annealing processor could not obtain a near-optimal solution (33,337) \cite{best_K2000}, the proposed method could obtain the near-optimal solution.

\begin{table}[t]
	\caption{Comparison of annealing accuracy with annealing processors on K2000.}
	\label{tb:annealing processors}
	\centering
	\begin{tabular}{|c||c|c|c|c|}
		\hline
		& SC-SA & CIM \cite{CIM} & SB \cite{SB} & STATICA \cite{STATICA} \\
		\hhline{|=#=|=|=|=|}
		Avg. & 33,262 & 32,459 & 32,768 & 33,073\\
		\hline
		Max  & 33,337 & 33,191 & N/ A & N/ A \\
		\hline
	\end{tabular}
\end{table}

\section{Conclusion}
\label{sec:conclusion}
The proposed method approximates the p-bits (spins) using stochastic computing, which can search for solutions of the combinatorial optimization problems around the global minimum energy.
As a result, the proposed SA method achieves the convergence speed a few orders of magnitude higher than conventional SA in combinatorial optimization problems: the MAX-CUT problems. 
Compared with conventional SA, the proposed SC-SA converged at 1,000 times the number of cycles and is 650 times faster in the MAX-CUT problem (K2000).
Compared with existing annealing processors, the proposed method achieves the best average cut value in the MAX-CUT problem (K2000).
Prospects for future research include simulating other combinatorial optimization problems using the proposed SC-SA method to evaluate the effectiveness of the proposed method. A large-scale hardware implementation of the proposed annealing method would be interesting as a fast solver of real-world social issues represented as the combinatorial optimization problems.

\section*{Acknowledgment}
This work was supported in part by JST CREST Grant Number JPMJCR19K3 and JSPS KAKENHI Grant Number JP21H03404

\bibliographystyle{IEEEtran}
\bibliography{sankou}

\end{document}